\documentstyle{amsppt}
\NoBlackBoxes
\magnification1200
\pagewidth{6.5 true in}
\pageheight{9 true in}

\topmatter
\title
Strong multiplicity one for the Selberg class
\endtitle
\author
K.  Soundararajan
\endauthor
\thanks{The author is partially supported by the American Institute of
Mathematics and by the
National Science Foundation.}\endthanks
\address
Department of Mathematics,  University of Michigan,
Ann Arbor, MI 48109, USA
\endaddress
\email
ksound\@umich.edu
\endemail
\endtopmatter

\def\phi{\varphi}
\def\lam{\lambda}

\document
\noindent In [7] A. Selberg axiomatized properties expected of $L$-functions
and introduced the ``Selberg class'' which is expected to coincide with
the class of all arithmetically interesting $L$-functions.   We recall
that an element $F$ of the Selberg class ${\Cal S}$ satisfies the
following axioms.

\noindent $\bullet$ In the half-plane $\sigma>1$ the
function $F(s)$ is given by a Dirichlet
series $\sum_{n=1}^{\infty} a_F(n) n^{-s}$ with $a_F(1)=1$ and
$a_F(n) \ll_{\epsilon} n^{\epsilon}$ for every $\epsilon >0$.

\noindent $\bullet$ There is a natural number $m_F$
such that $(s-1)^{m_F} F(s)$ extends to an analytic function in the
entire complex plane.

\noindent $\bullet$  There is a function $\Phi_{F}(s) = Q_F^{s} \Gamma_F(s)
F(s)$ where $Q_F >0$ and
$$
\Gamma_F(s) = \prod_{j=1}^{r_F} \Gamma(\lam_j(F)s+\mu_j(F))
\qquad \text{with} \qquad \lam_j(F)>0 \text{  and  } \text{Re }\mu_j(F) \ge 0
$$
such that
$$
\Phi_F(s) = \omega_F \overline{\Phi}_F(1-s),
$$
where $|\omega_F| =1$ and for any function $f$ we denote $\overline{f}(s)
= \overline{f(\overline{s})}$. We let $d_F:= 2\sum_{j=1}^{r} \lam_j$ denote
the ``degree'' of $F$.

\noindent $\bullet$ We may express $\log F(s)$ by a Dirichlet series
$\log F(s) = \sum_{n=2}^{\infty} b_F(n) \Lambda(n)/(n^s \log n)$
where $b_F(n) \ll  n^{\vartheta}$ for some $\vartheta <\frac 12$.  We
adopt the convention that $b_F(n)=0$ if $n$ is not a prime power.  From
the assumption $a_F(n) \ll_{\epsilon} n^{\epsilon}$ it follows that
$b_F(p^k) \ll_{k,\epsilon} p^{\epsilon}$.

It is believed that the Selberg class satisfies the
following ``strong multiplicity one'' principle:
If $F$ and $G$ are two elements of the
Selberg class with $a_F(p)=a_G(p)$ (equivalently $b_F(p)=b_G(p)$)
for all but finitely many primes $p$ then $F=G$.
In [5] R. Murty and K. Murty prove such a result under the
additional hypothesis that $a_F(p^2)=a_G(p^2)$ for all but finitely
many primes $p$.  Recently J. Kaczorowski and A. Perelli [2] have
established this principle under the additional hypothesis
that
$$
\lim_{\sigma \to 1^+} (\sigma-1) \sum_p \frac{|a_F(p^2)-a_G(p^2)|}{p^{\sigma}}
\log p < \infty.
$$
This criterion is equivalent to saying that $a_F(p^2) - a_G(p^2)$ is
bounded on average: that is for all large $x$
$$
\sum_{p \le e^{x}} \frac{|a_F(p^2)- a_G(p^2)|^2}{p} \log p \ll x. \tag{1}
$$

In this note we develop a different method which establishes
the strong multiplicity one principle under a much weaker hypothesis
than (1), but which fails (barely) to prove the full principle.
We use $\log_j$ to denote the $j$-fold iterated logarithm; thus
$\log_2 =\log \log$, $\log_3 =\log \log \log $ and so on.

\proclaim{Theorem}  Suppose $F$ and $G$ are elements of the
Selberg class with $a_F(p)=a_G(p)$ (equivalently, $b_F(p)=b_G(p)$)
for all $p \notin {\Cal E}$ where ${\Cal E}$ is a thin
set of primes satisfying
$$
 \# \{p \in {\Cal E}, p\le x\} \ll x^{\frac 12 -\delta} \tag{2}
$$
for some fixed $\delta >0$.  Then $F$ and $G$ have
the same degree:  $d_F=d_G$.  If in addition we have
$$
\sum_{p\le e^x} \frac{|a_F(p^2)-a_G(p^2)|^2}{p} \log p \ll \exp\Big(
\frac{x}{\log x (\log_2 x)^5}\Big),
\tag{3}
$$
then $F=G$.
\endproclaim

Although (3) is a considerably weaker restriction than (1)
it is still stronger than the bound $\ll e^{\epsilon x}$ which is all we
know in general.  When combined with the classification of
elements of the Selberg class of degree 1 by Kaczorowski and
Perelli (see [3]) our Theorem permits the following corollary.

\proclaim{Corollary}  Suppose $\chi$ is a primitive Dirichlet
character and that $F$ is an element of the Selberg
class with $f(p)=\chi(p)$ for all $p\notin {\Cal E}$ where
${\Cal E}$ is a set of primes satisfying (2).  Then $F(s)=L(s,\chi)$.
\endproclaim

To deduce the Corollary we note that our Theorem implies that the
degree of $F$ is $1$.  Since Kaczorowski and Perelli [3] have
shown that the only elements of the Selberg class
of degree $1$ are Dirichlet $L$-functions the Corollary follows.

We now embark on the proof of our Theorem.
Put $c(n) =b_F(n)-b_G(n)$ and observe that in $\text{Re } s> \frac 32$ we
have
$$
-\frac{F^{\prime}}{F}(s) + \frac{G^{\prime}}{G}(s) = \sum_{n=1}^{\infty}
\frac{c(n) \Lambda(n)}{n^s}
= \sum_{k=1}^{\infty} \sum_{p} \frac{c(p^k) \log p}{p^{ks}}. \tag{4}
$$
Since $c(p)=0$ unless $p \in {\Cal E}$ we see from (2) that
$\sum_p c(p)\log p/p^s$ is entire in Re$(s) > \frac 12-\delta$.
Since $c(p^2)$ and $c(p^3)$ are $\ll p^\epsilon$ we see that
$\sum_p c(p^2)\log p/p^{2s}$ and $\sum_p c(p^3)\log p/p^{3s}$ are
entire in $\text{Re } s >\frac 12$ and $\frac 13$ respectively.  Lastly
since $c(n) \ll n^{\vartheta}$ for $\vartheta <\frac 12$ it
follows that $\sum_{k=4}^{\infty} \sum_p c(p^k)\log p/p^{ks}$ is
entire in $\text{Re }s >\frac 12$.  Thus the RHS of (4) is an
entire function in Re $s > \frac 12$ furnishing an analytic
continuation of $-\frac{F^{\prime}}{F}(s) +\frac{G^{\prime}}{G}(s)$ in
this region.  It follows that the zeros of $F$ and $G$ in this
region coincide (including multiplicities), and also that $F$ and
$G$ have poles of the same order at $1$ (that is, $m_F=m_G$).
Arguing similarly with $\overline{F}(s)=\sum_{n} \overline{a_F(n)}{n^{-s}}$
and  $\overline{G}(s)$ we see that their zeros in Re $s>\frac 12$ also
coincide.  Using the functional equation it follows that
the zeros of $\Phi_F(s)$ and $\Phi_G(s)$ in Re $s<\frac 12$ also coincide.
Summarizing we see that $\Phi_F(s)$ and $\Phi_G(s)$ have the same
zeros except possibly on the critical line Re $s=\frac 12$, and that
they have (possibly) poles of the same order at $s=1$.

Let $\rho_F =\frac 12+ i\gamma_F$ and $\rho_G =\frac 12+i\gamma_G$
denote typical zeros of $\Phi_F(s)$ and $\Phi_G(s)$.  We do not
suppose that $\gamma_F$ and $\gamma_G$ are real, although this version
of the Riemann hypothesis is expected to be true.
A standard application of the argument principle shows that
 $\#\{\rho_F : \ \ |\text{Im }\rho_F|\le T\} = \frac{d_F}{\pi} T\log T
+ c_F T +O_F(\log T)$ where $d_F$ is the degree and $c_F$ is a constant.
Similar estimates apply for the zeros $\rho_G$ up to height $T$.

We now recall an explicit formula connecting the zeros $\rho_F$ to
the prime power values $b_F(n)\Lambda(n)$; for details see
for example Z. Rudnick and P. Sarnak [6].  Let $g$ be a smooth
compactly supported function and put $h(s) =\int_{-\infty}^{\infty} g(u)
e^{isu} du$.  We may recover $g$ from $h$ by means of
the Fourier inversion formula $g(x) =\frac{1}{2\pi }\int_{-\infty}^{\infty}
 h(u)e^{-ixu} du$.  The explicit formula now reads
$$
\align
\sum_{\gamma_F} h(\gamma_F)
&= m_F \Big( h\Big(-\frac i2\Big) +h\Big(\frac i2\Big)\Big)
\\
&+ \frac{1}{2\pi} \int_{-\infty}^{\infty} h(r)
\Big( 2\log Q_F +\frac{\Gamma_F^{\prime}}{\Gamma_F}(\tfrac 12 +ir)
+ \frac{\overline{\Gamma_F}^{\prime}}{\overline{\Gamma_F}} (\tfrac 12-ir) \Big)
dr
\\
&- \sum_{n=1}^{\infty} \Big( \frac{b_F(n)\Lambda(n)}{\sqrt{n}} g(\log n)
+ \frac{\overline{b_F(n)}\Lambda(n)}{\sqrt{n}} g(-\log n)\Big).
\tag{5}
\\
\endalign
$$

We now choose $g$ to be a smooth non-negative function in $[-1,1]$ such that
$g(0)=\frac{1}{2\pi }\int_{-\infty}^{\infty} h(r)dr \gg 1$ and
such that
$$
|h(t)| \ll \exp(- |t|/\log^2 |t|), \tag{6}
$$
for large $|t|$.  A result of A.E. Ingham [1] shows that such a choice of
$g$ is possible.  For example, following Ingham we may take
$$
h(t) = \prod_{n=N}^{\infty} \Big(\frac{\sin(2\pi t/(n (\log n)^{\frac 32}))}{
2\pi t/(n (\log n)^{\frac 32})}\Big)^2
$$
for some large $N$ and then one could check that $h$ and its
Fourier transform $g$ are a permissible choice.


Let $T$ be a large positive number and let $t$ be in $(T,2T)$.
Let $\log^2 T \ge L\ge \log T$
be a large parameter to be chosen later.
Then an application of (5) gives that
$$
\align
\sum_{\gamma_F} h(L(\gamma_F-t))
&= m_F \Big( h(L(-\tfrac i2-t)) +
h(L(\tfrac i2+t))\Big) \\
&+ \frac{1}{2\pi}
\int_{-\infty}^{\infty} h(L(r-t))
\Big( 2\log Q_F +\frac{\Gamma_F^{\prime}}{\Gamma_F}(\tfrac 12 +ir)
+ \frac{\overline{\Gamma_F}^{\prime}}{\overline{\Gamma_F}} (\tfrac 12-ir) \Big)
dr
\\
&- \frac{1}{L}
\sum_{n=1}^{\infty} \Big( \frac{b_F(n)\Lambda(n)}{n^{\frac 12+it}}
g\Big(\frac{\log n}{L}\Big)
+ \frac{\overline{b_F(n)}\Lambda(n)}{n^{\frac 12-it}} g\Big(-\frac{\log n}{L}
\Big)\Big).
\tag{7}
\\
\endalign
$$
We call the middle term on the RHS above $H_F(t,L)$ and the third
term there $D_F(t,L)$.  Stirling's
formula shows that
$$
H_F(t,L) = \frac{d_F \log T +O(1)}{L}
\frac{1}{2\pi} \int_{-\infty}^{\infty} h(r) dr
= g(0) \frac{d_{F}\log T}{L} + O(1/L).
\tag{8}
$$

We apply the explicit formula (7) for $G$ as well and subtract the
two formulae.  From our remarks on the poles of $F$ and $G$
and their zeros not on the critical line we conclude that
$$
Z_F(t,L) -Z_G(t,L) = H_F(t,L)-H_G(t,L) - D_F(t,L)+D_G(t,L),
\tag{9a}
$$
where
$$
Z_F(t,L) =\sum_{\gamma_F \in {\Bbb R}} h(L(\gamma_F -t))
\qquad \text{and } \qquad Z_G(t,L) =\sum_{\gamma_G \in {\Bbb R}}
h(L(\gamma_G-t)). \tag{9b}
$$

We now record a mean-value estimate for the $D_F(t,L)-D_G(t,L)$ terms.
Since $|c(p)|\ll p^{\epsilon}$ we
see from (2) that $\sum_{p} |c(p)g(\log p/L)|\log p/\sqrt{p}
= O(1)$.  Further as noted earlier $\sum_{k\ge 3} \sum_p |c(p^k)\log p\,
g(k\log p/L)|/p^{k/2} = O(1)$.  Thus
$$
\align
L (D_F(t,L)-D_G(t,L))
&\ll 1 +
\Big|
\sum_p \frac{a_F(p^2) -a_G(p^2)}{p^{1+2it}}
\log p \, g\Big(\frac{2\log p}{L}\Big)\Big|
\\
&\hskip .5 in +\Big|
\sum_p \frac{\overline{a_F(p^2)} -\overline{a_G(p^2)}}{p^{1-2it}}
\log p \, g\Big(-\frac{2\log p}{L}\Big)\Big|.
\\
\endalign
$$
Using now a familiar mean-value estimate of Montgomery and Vaughan ([4], see
Corollary 3) we see that
$$
\align
\int_T^{2T}
\Big| \sum_p \frac{a_F(p^2) -a_G(p^2)}{p^{1+2it}}
\log p \, g\Big(\frac{2\log p}{L}\Big)\Big|^2 dt
&\ll \sum_{p\le e^{L/2}} \frac{|a_F(p^2)-a_G(p^2)|^2}{p^2} (T+p)\log^2 p
\\
&\ll T+ \sum_{p\le e^{L/2}} \frac{|a_F(p^2)-a_G(p^2)|^2}{p}\log^2 p.
\\
\endalign
$$
We conclude unconditionally that
$$
\int_T^{2T} (L|D_F(t,L)-D_G(t,L)|)^2 dt \ll T+ e^{L\epsilon}, \tag{10a}
$$
and if we assume the condition (3) that
$$
\int_T^{2T} (L|D_F(t,L)-D_G(t,L)|)^2 dt \ll T+ \exp\Big(\frac{L}
{\log L (\log_2 L)^5}\Big). \tag{10b}
$$

Let $W \ge 1$ be a real parameter and let ${\Cal L}={\Cal L}(W)$ denote
the set of $t\in [T,2T]$ such that there exists either $\gamma_F$ or
$\gamma_G$ in $(t-1/(W\log T), t+1/(W\log T))$.  Let $\overline{{\Cal L}}
$ denote the complementary set $[T,2T]\backslash{\Cal L}$.  Since
there are $\ll T\log T$ ordinates $\gamma_F$ or $\gamma_G$ in $[T,2T]$
we see that $\text{meas} ({\Cal L}) \ll T/W$.  Now
$$
L \int_{t\in \overline{{\Cal L}}} |Z_F(t,L)| dt \ll L \int_T^{2T}
\sum\Sb \gamma_F \in {\Bbb R}\\ |\gamma_F-t|\ge 1/(W\log T)\endSb
|h(L(\gamma_F-t))| dt.
$$
If the distance of $\gamma_F$ from $(T,2T)$ exceeds $n$ then
by (6) the contribution of such a $\gamma_F$ to the RHS above
is $\ll  \exp(-nL/(\log nL)^2)$.  Further the number of
ordinates $\gamma_F$ whose distance from $(T,2T)$ is between $n$ and $n+1$
is $\ll \log (T(n+1))$ and so we conclude that the contribution to
the RHS above from zeros not in $(T-1,2T+1)$ is $\ll 1$. Thus
$$
\align
L \int_{t\in \overline{{\Cal L}}} |Z_F(t,L)| dt
&\ll 1+\sum_{\gamma_F \in (T-1,2T+1)} \int_{|y|\ge L/(W\log T)} |h(y)|dy
\\
&\ll 1+ T\log T \int_{|y|\ge L/(W\log T)} |h(y)|dy.
\\
\endalign
$$
A similar estimate applies for $|Z_G(t,L)|$ so that
$$
L \int_{t\in \overline{{\Cal L}}} |Z_F(t,L) -Z_G(t,L)| dt \ll 1+
 T\log T \int_{|y|\ge L/(W\log T)} |h(y)|dy.  \tag{11}
$$

With these preliminaries in place we are now ready to finish the
proof of our Theorem. Suppose first that $d_F \neq d_G$.
Then by (8) we know that $L(H_F(t,L)-H_G(t,L)) \gg \log T$.  Integrating
(9a) over $t\in \overline{{\Cal L}}$ we find that
$$
\align
T\log T &\ll \int_{t\in \overline{{\Cal L}}} L|H_F((t,L)-H_G(t,L)| dt
\\
&\ll \int_{t \in \overline{{\Cal L}}} L|Z_F(t,L)-Z_G(t,L)| dt +
\int_{t \in \overline{{\Cal L}}} L|D_F(t,L)-D_G(t,L)| dt
\\
&\ll 1+ T\log T \int_{|y|\ge L/(W\log T)} |h(y)|dy + (T+e^{L\epsilon}),\\
\endalign
$$
using (11) and (10a) above.  We now choose $L=W^{2} \log T$
and choose $W$ to be large but smaller than $1/\sqrt{\epsilon}$.  Then for
large $T$ the above gives that $1 \ll \int_{|y|>W} |h(y)| dy$
which is impossible for sufficiently large $W$ in view of (6).

We may now suppose that $d_F =d_G$.  In this case we will require
the additional hypothesis (3) which permits the bound (10b) above.  Let
$m$ be a fixed integer with $c(m) \neq 0$.  Naturally if
no such $m$ exists then $F=G$.  We multiply (9a) by $Lm^{it}$ and integrate
for $t\in \overline{{\Cal L}}$.  Using (11) we see that the LHS gives
$$
\int_{t \in \overline{{\Cal L}}} Lm^{it} (Z_F(t,L)-Z_G(t,L)) dt
\ll 1+ T\log T \int_{|y| \ge L/(W\log T)} |h(y)| dy.  \tag{12a}
$$
On the other hand using the RHS of (9a) this is also equal to
$$
\int_{t \in \overline{{\Cal L}}} Lm^{it} (H_F(t,L)-H_G(t,L)-D_F(t,L) +
D_G(t,L)) dt.
$$
Now
$$
\align
&\int_{t \in \overline{{\Cal L}}} Lm^{it} (H_F(t,L)-H_G(t,L)) dt
\\
&\hskip .3 in
= \int_{T}^{2T} Lm^{it} (H_F(t,L)-H_G(t,L)) dt +O\Big( \int_{t\in{\Cal L}}
L|H_F(t,L)-H_G(t,L)| dt \Big). \\
\endalign
$$
Since $d_F=d_G$ we see by (8) that the second term above is
$\ll \text{meas}({\Cal L}) \ll T/W$.  Further using integration
by parts, and since $L\frac{d}{dt} (H_F(t,L)-H_G(t,L)) \ll 1/T$
by Stirling's formula, we see that the first term above is $\ll 1$.
Thus
$$
\int_{t \in \overline{{\Cal L}}} Lm^{it} (H_F(t,L)-H_G(t,L)) dt \ll
1+T/W. \tag{12b}
$$
Further
$$
\align
&\int_{t \in \overline{{\Cal L}}} Lm^{it} (D_F(t,L)-D_G(t,L)) dt
\\
&\hskip .75 in = \int_{T}^{2T} Lm^{it} (D_F(t,L)-D_G(t,L)) dt + O\Big(
\int_{t\in {\Cal L}} L|D_F(t,L)-D_G(t,L)|dt\Big).
\\
\endalign
$$
The second term above is $\ll (T+\exp(\frac{L}{\log L(\log_2 L)^5})
/\sqrt{W}$ upon using
(10b) and Cauchy's inequality.  Integrating term by term we see that
the first term is
$$
T \frac{c(m)}{\sqrt{m}} \Lambda(m) + O\Big(\sum_{n\le e^{L}}
\frac{|c(n)|\log n}{\sqrt{n}} \Big) = T \frac{c(m)}{\sqrt{m}} \Lambda(m)
+ O\Big(\exp\Big(\frac{L}{\log L (\log_2 L)^5}\Big)\Big)
$$
using (2) and (3).  Thus we conclude that
$$
\int_{t \in \overline{{\Cal L}}} Lm^{it} (D_F(t,L)-D_G(t,L)) dt
= T \frac{c(m)}{\sqrt{m}} \Lambda(m) + O\Big(\frac{T}{\sqrt{W}} +
\exp\Big(\frac{L}{\log L (\log_2 L)^5}\Big)
\Big). \tag{12c}
$$
Combining (12a,b,c) we get that
$$
1\ll \frac{c(m)}{\sqrt{m}}\Lambda(m) \ll \frac{1}{T} + \log T
\int_{|y| \ge L/(W\log T)} |h(y)| dy + \frac{1}{\sqrt{W}} +
\frac{\exp(\frac{L}{\log L (\log _2 L)^5})}{T}.
$$
We now choose $W=\log_3 T$, $L=\log T \log_2 T (\log_3 T)^4$ and use (6) to
obtain a contradiction.

{\bf Acknowledgements.}  I am grateful to Hugh Montgomery for referring
me to Ingham's paper [1].


\Refs

\widestnumber \key{1}

\ref\no 1
\by A.E. Ingham
\paper A note on Fourier transforms
\jour J. London Math. Soc.
\vol 9
\yr 1934
\pages 29--32
\endref

\ref\no 2
\by J. Kaczorowski and A. Perelli
\paper Strong multiplicity one for the Selberg class
\jour  Comptes Rendus Acad. Sci. Paris Ser. I Math.
\vol 332
\yr 2001
\pages 963--968
\endref

\ref \no 3
\by J. Kaczorowski and A. Perelli
\paper On the structure of the Selberg class. I. $0\le d\le 1$
\jour Acta Math.
\vol 182
\yr 1999
\pages 207--241
\endref


\ref
\no 4
\by H.L. Montgomery and R.C. Vaughan
\paper Hilbert's inequality
\jour J. Lond. Math. Soc. (2)
\vol 8
\yr 1974
\pages 73--82
\endref

\ref\no 5
\by M.R. Murty and V.K. Murty
\paper Strong multiplicity one for the Selberg class
\jour Comptes Rendus Acad. Sci. Paris Ser. I Math.
\vol 319
\yr 1994
\pages 315--320
\endref

\ref\no 6
\by Z. Rudnick and P. Sarnak
\paper Zeros of principal $L$-functions and random matrix theory
\jour Duke Math. J.
\vol 81
\yr 1996
\pages 269--322
\endref

\ref \no 7
\by A. Selberg
\paper Old and new conjectures about a class of Dirichlet series
\jour Collected papers
\vol II
\pages 47--63
\publ Springer, Berlin (1991)
\endref


\endRefs
\enddocument